\documentclass[14pt,a4paper]{article}

\usepackage[utf8]{inputenc}
\usepackage[T2A]{fontenc}
\usepackage[english,russian]{babel}

\usepackage{indentfirst} % Красная строка
\usepackage[14pt]{extsizes}
\usepackage{amssymb}
\usepackage{amsmath}
\usepackage{amsthm}
\usepackage{pifont}
\usepackage{graphicx}
\usepackage[colorinlistoftodos]{todonotes}
\usepackage[colorlinks=true, allcolors=blue]{hyperref}
\usepackage{algorithm}
\usepackage{algorithmic}

\newtheorem{theorem}{Теорема}

\def\UDK#1{{\leftline{УДК {#1}}}}

\usepackage{geometry}
\newcommand{\demyan}[1]{\textcolor{purple}{#1}}

\begin{document}
\renewcommand{\abstractname}{\vspace{-\baselineskip}}

$$
\\\\
$$

\UDK{519.853.62}

\begin{center}
\textbf{О многостадийной модели равновесного распределения транспортных потоков по путям и достаточных условиях, когда поиск равновесия сводится к решению задачи оптимизации}

Е.\,В.~Гасникова, А.\,В.~Гасников, Д.\,В.~Ярмошик, М.\,Б.~Кубентаева, М.\,И.~Персиянов, И.\,В.~Подлипнова, Е.\,В.~ Котлярова, И.\,А.~Склонин, Е.\,Д.~Подобная, В.\,В.~Матюхин

% \textbf{Д.\,Камзолов, А.\,И.~Тюрин,  Д.\,А.~Пасечнюк}

\end{center}

\begin{abstract}
\noindent Многостадийное моделирование транспортных потоков начало активно развиваться с 70-х годов прошлого века. Были созданы пакеты транспортного моделирования, в основе которых имеется набор задач выпуклой оптимизации, последовательное отрешивание которых (по циклу) приближает к искомому равновесному распределению. Альтернативный путь -- попробовать найти такую общую задачу выпуклой оптимизации, решение которой давало бы искомое равновесие. В отличие от первого пути, по второму пути не всегда удается пройти. В данной статье предпринята попытка найти достаточные условия, гарантирующие, что второй путь будет успешен. Общность, которая выбрана в данной статье, насколько нам известно, ранее не встречалась в литературе.
В частности, возможность в качестве одного из блоков многостадийной модели выбирать модель стабильной динамики (а не общепринятую модель Бэкмана) в сочетании с возможностью выбирать различные типы пользователей и транспортных средств является новой. %тут кажется не нужна запятая
\end{abstract}

\section{Введение} \label{section_1}
В данной работе описывается  вариационный (экстремальный) принцип, сводящий при определенных условиях поиск равновесия в многостадийной модели распределения транспортных потоков по сети к задаче выпуклой оптимизации. Под распределением потоков, прежде всего, понимается: 1) расчёт матрицы корреспонденций и 2) распределение потоков по путям при заданных корреспонденциях. Таким образом, речь идет о двухуровневой модели распределения. Многостадийные модели транспортных потоков являются одним из основных объектов изучения при долгосрочном транспортном планировании \cite{Ortuzar2002,gasnikov2020posobie,Boyles2020}. С помощью таких моделей можно просчитывать долгосрочные последствия: введения в эксплуатацию различных инфраструктурных объектов, изменения дорожной сети и т.п.
В отличие от подхода \cite{kotlyarova2021} здесь мы погружаем в двухстадийную модель полноценную четырехстадийную модель, которая, среди прочего, учитывает расщепление перемещений на личном и общественном транспорте, а также может учитывать множество других транспортных средств и сетей дорог. Обсуждаются различные слои спроса. 

Следуя \cite{gasnikov2014matmod}, в статье выписывается выпукло-вогнутая седловая задача, к которой сводится поиск равновесия в такой двухстадийной модели. Далее эта задача упрощается (путем перехода к двойственному представлению), и совсем кратко описывается численный способ решения возникающей в итоге (двойственной) задачи. Отличительными особенностями данной работы являются: простой способ получения итоговой задачи оптимизации (седловой задачи) и способ ее решения. В предлагаемых ранее подходах либо не рассматривалась возможность предельного перехода от модели Бэкмана к модели стабильной динамики (Нестеров--деПальма) \cite{Evans1976,DeCea2005}, либо рассматривался только один слой спроса и одна транспортная сеть \cite{gasnikov2014matmod}. 

\section{Основные определения и обозначения}\label{section_2}
Для простоты будем рассматривать замкнутую транспортную систему, описываемую графом $G = \langle V, E\rangle$, где $V$ -- множество вершин ($|V| = n$), а $E$~-- множество ребер ($|E| = m$). Будем обозначать ребра графа через $e\in E$. Для стандартной транспортной системы можно ожидать, что $m\simeq 3n$. Для больших мегаполисов (таких как Москва)  $n \simeq 10^5$. 
Транспортный граф $G$ считается известным. 

Будем считать, что имеются разные слои спроса, которые мы будем индексировать буквой $r$. Например, слой спроса, отвечающий передвижениям <<дом-работа>> или <<дом-учеба>>. Внутри каждого слоя спроса допустимы разные типы пользователей $t \in M(r)$ ($M(r)$ -- множество типов пользователей, отвечающих слою спроса $r$), например, имеющие личный транспорт, которые могут выбирать между личным и общественным транспортом и класс пользователей, которые не имеют личный транспорт. Тип транспорта будем обозначать индексом $k$. Множество типов транспортных средств, доступных пользователям типа $k$ будем обозначать $Z(t)$. Множество всех типов транспортных средств $K$ можно разбить на классы транспортных средств $\left\{K_{b}\right\}$, использующих (внутри класса) одну и ту же транспортную сеть $b$ ($b\in B$). Под транспортной сетью понимается набор маршрутов $P^b$ (составленных из ребер $G$), доступных для перемещения по данной сети. Например, можно выделить обычную сеть дорог для личного транспорта, которой пользуются легковые автомобили и частично общественный транспорт. Можно выделить велодорожки (особенно популярно в Германии, Бельгии и Голландии), которыми могут пользоваться не только велосипедисты, но, например, и самокатные средства.

Часть вершин $O\subseteq V$ (\textit{origin}) являются источниками корреспонденций, а часть \demyan{--} стоками корреспонденций $D\subseteq V$ (\textit{destination}). Если говорить более точно, то вводится множество пар (источник, сток) корреспонденций $OD \subseteq V\bigotimes V$. Сами корреспонденции будем обозначать через $d^{rt}_{ij}$, где $(i,j)\in OD$. Как правило $|OD|\ll n^2$ \cite{gasnikov2014matmod}. Не ограничивая общности, будем далее считать, что $\sum_{r,t\in M(r);(i,j)\in OD} d^{rt}_{ij} = 1$. Множество пар $OD$ считается известным. Корреспонденции -- не известны! Однако известны (заданы) характеристики источников и стоков корреспонденций, соответственно, \demyan{$l_i^{r}$ и $w_j$}. То есть известны величины \cite{Wilson1978} \demyan{$\{l_i^{r}\}_{i\in O}$, $\{w_j\}_{j\in D}$}:
\begin{equation}\label{corr}
    % \sum_{k \in Z(t)}
    \demyan{\sum_{t\in M(r)}} \sum_{j: (i,j)\in OD} d^{rt}_{ij} = l_i^{\demyan{r}}, \quad \demyan{\sum_{r}}\sum_{t\in M(r)}\sum_{i: (i,j)\in OD} d^{rt}_{ij} = \demyan{w_j}.
\end{equation}
Заметим, что \demyan{$\sum_{r}\sum_{i\in O} l^{r}_i = \sum_{j\in D} w_j = 1$}. Условие \eqref{corr} будем также для краткости записывать в виде $d\in (l,w)$.

Обозначим через (зависимость $b(k)$ определяет тип сети по типу транспортного средства: $b = b(k)$ тогда и только тогда, когда $k \in K_b$) $$\tau_e^k(f_{e,b(k)},f_{e}^k)=\tau_e(f_{e,b(k)}) + \tilde{\tau}^k_e(f_{e}^k)$$  функцию затрат (например, временных) на проезд по ребру (участку дороги) $e$ на типе транспортных средств $k$ на сети $b(k)$, если суммарный поток  на этом участке сети $f_{e,b(k)}$, а поток транспортных средств типа $k$ равен $f_{e}^k$ ($f_{e,b(k)} = \sum_{k \in K_b} f_{e}^k$). Функции $\tau_e(f_{e,b(k)})$ и $\tilde{\tau}^k_e(f_{e}^k)$ считаются заданными, например, таким образом: \cite{gasnikov2013book,Patriksson2015}
\begin{equation}\label{BPR}
    \tau_e(f) = \bar{t}_e\left(1 +\zeta\left(\frac{f}{\bar{f}_e}\right)\right)^{\frac{1}{\mu}}, 
\end{equation}
где $\bar{t}_e$ -- время прохождения ребра $e$, когда участок свободный (определяется разрешенной скоростью на данном участке), а $\bar{f}_e$ -- пропускная способность ребра $e$ (определяется полосностью: [пропускная способность] $\le$ [число полос] * [2000 авт/час] и характеристиками перекрестков).
Считается, что эти характеристики известны \cite{Stabler2020}. Параметр $\mu = 0.25$ -- BPR-функции \cite{Patriksson2015}, но допускается и $\mu\to 0+$ -- модель стабильной динамики \cite{Nesterov2003,gasnikov2014matmod,kotlyarova2022}. Параметр $\zeta >0$ также считается заданным. 
Относительно $\tilde{\tau}^k_e(f)$  часто считают, что $\tilde{\tau}^k_e(f) \equiv c_e^k$ \cite{Shvetsov2003}.
% \todo[inline]{Мб $\tilde{\tau}^k_e(f) = c_e^k$?  Вроде $c_e^k f $~-- это $\sigma$, a не $\tau$.\\
% Ниже ещё написано ``Отметим также, что если $\tilde{\tau}^k_e(f) = c_e^k f$, то можно считать $\tilde{t}_e^k = c_e^k$''. Получается странно, разве не должно быть $\tau =t$?
% }
Например, если транспортное средства $k$ весит больше $3.5$ тонн, то по некоторым ребрам $e$ проезд может быть запрещен и тогда можно полагать $c_e^k=\infty$.

Полезно также ввести $t_e^k$ -- (временные) затраты на прохождения ребра $e$ на транспортном средстве типа $k$. Обозначение $t$ ранее уже было использовано под тип пользователя, однако это не должно вызывать недоразумение, поскольку далее по смыслу будет понятно, какое $t$ имеется в виду в том или ином случае.  Согласно вышенаписанному 
\begin{equation}\label{t(f)}
   t_e^k = t_{e,b(k)}+\tilde{t}_{e}^k=\tau_e(f_{e,b(k)}) + \tilde{\tau}^k_e(f_{e}^k)=\tau_e^k(f_{e,b(k)},f_{e}^k). 
\end{equation}
По этим затратам $t^k = \{t_e^k\}_{e\in E}$ можно определить затраты на перемещение из источника $i$ в сток $j$ по кратчайшему пути: 
\begin{equation}\label{T}
   T^k_{ij}(t^k) = \min_{p \in P^{b(k)}_{ij}} \sum_{e\in E} \delta_{ep}t_e^k, 
\end{equation}
 где $p$ -- путь (без самопересечений -- циклов) на графе (набор ребер),
$P_{ij}^{b(k)}$~-- множество всевозможных путей на графе транспортной сети $b(k)$, стартующих из источника $i$ и заканчивающихся в стоке $j$, $\delta_{ep} = 1$, если ребро $e$ принадлежит пути $p$ и $\delta_{ep} = 0$ -- иначе.
% В зависимости от контекста, в эти обозначения будут добавляться верхние индексы, отвечающие типу пассажиров $k$.

В ряде выкладок далее также будет полезен вектор $x^{rtk} = \{x_p^{rtk}\}_{p \in P^{b(k)}}$~-- вектор распределения потоков по путям, где $P^b = \bigcup_{(i,j)\in OD} P_{ij}^b$. Заметим, что 
\begin{equation}\label{f(x)}
 f_e^k = \sum_{r;~t\in M(r);~p\in P^b(k)} \delta_{ep} x_p^{rtk}.
\end{equation}
% или в матричном виде $f^k = \Theta x^{rtk}$, где $\Theta = \|\delta_{ep}\|_{e\in E, p\in P}$.

% Также определим $T^{tk}_{ij}$ как средние (временные) затраты для пассажиров типа $t$, использующих тип транспорта $k$ (не важно какого именно слоя спроса), которые требуются, чтобы добраться из района $i$ в район $j$. Впоследствии мы свяжем эти затраты с введенными ранее функциями $T^{tk}_{ij}(t)$. Если $|Z(t)| = 1$, то будем писать $T^{t*}_{ij}(t)$.

Все дальнейшие обозначения будем вводить по мере надобности.

\section{Энтропийная модель А.Дж.~Вильсона расчёта матрицы корреспонденицй}

Под энтропийной моделью расчета матриц корреспонденций $d^{rt}(\bar{T})$ понимается определенный способ вычисления набора корреспонденций $\{d^{rt}_{ij}\}_{(i,j)\in OD}$ по известным матрицам средних затрат $\bar{T} = \{\bar{T}^{t}_{ij}\}_{(i,j)\in OD}$. Этот способ заключается в решении задачи энтропийно-линейного программирования, которую можно понимать, как энтропийно-регуляризованную транспортную задачу
% \footnote{Вместо $d^{rt}_{ij}\ln d^{rt}_{ij}$  точнее было бы писать $d^{rt}_{ij}\ln\left(d^{rt}_{ij}/\left(\sum_{r;~t\in M(r);~(i,j)\in OD}d^{rt}_{ij}\right)\right)$ \cite{gasnikov2014matmod}, но не ограничивая общности нормировку можно выбрать так, что $\sum_{r;~t\in M(r);~k\in Z(t);~(i,j)\in OD} d_{ij}^{rtk} = 1$, поэтому, возможно, упрощенная форма записи.}
\begin{equation}\label{Wilson}
   \min_{d\in (l,w);d\ge 0} \left\{\sum_{r}\beta^r\sum_{t\in M(r);~(i,j)\in OD} d^{rt}_{ij}\bar{T}_{ij}^{t} +  \sum_{r;~t\in M(r);~(i,j)\in OD} d^{rt}_{ij}\ln d^{rt}_{ij}\right\},
\end{equation}
где параметры $\beta^r >0$ считаются известными \cite{Wilson1978,gasnikov2013book,gasnikov2020posobie}. Относительно выбора этих параметров, см. \cite{gasnikov2014matmod,gasnikov2020posobie,ivanova2020}. В действительности, это, так называемые, структурные параметры, и от их подгонки существенно зависит качество модели. Если известны средние времена $C^r$ в пути для разных слоев спроса $r$, то находить $\beta^r$ можно из системы уравнений
$$\sum_{t\in M(r);~(i,j)\in OD} d^{rt}_{ij}(\beta)\bar{T}_{ij}^{tk}\simeq C^r.$$

\textbf{Пример.} В качестве простого примера практической ситуации можно взять $r = \{1,2\}$, где $r=1$ -- отвечает слою спроса <<дом-работа>>, а $r = 2$ -- отвечает слою спроса <<дом-учеба>>. При этом пассажиры типа $t=1$ могут пользоваться личным ($k=1$) и общественным транспортом ($k=2$), $t=2$~-- только общественным, а $t=3$ возникает только в слое спроса $r=2$ и отвечает пассажирам, использующих  предусмотренный спец. транспорт, например, школьный автобус ($k=3$). 
% Пример немного искусственный, поскольку едва ли можно рассчитывать на массовость такого типа передвижений, но нужен нам, чтобы продемонстрировать <<эффект размера транспортного средства>>. 

К сожалению, если $\beta_r$ не равны между собой, то не получается эволюционно проинтерпретировать модель расчета матрицы корреспонденций, как некоторое равновесие микросистемы в стиле \cite{gasnikov2013book,gasnikov2014matmod,gasnikov2016,gasnikov2020posobie,Gasnikova2023}. Как мы увидим в дальнейшем, в этом случае также не получается свести поиск равновесия в многостадийной модели к задаче оптимизации.

\iffalse
Полезно также отметить, что решение задачи \eqref{Wilson} можно представить в виде \cite{Wilson1978}
\begin{align*}
    d^{rtk}_{ij}=A^{rt}_{i} B_j^r Q^{rt}_i D_j^r \exp\left(-\beta^r T^{tk}_{ij}\right),\\
    A^{rt}_i = \left(\sum_{k\in Z(t);~j} B_j^r D_j^r\exp\left(-\beta^r T^{tk}_{ij}\right)\right)^{-1},\\
    B_j^r = \left(\sum_{t\in M(r);~k\in Z(t);~j} A^{rt}_i Q^{rt}_i\exp\left(-\beta^r T^{tk}_{ij}\right)\right)^{-1}.
\end{align*}
Отсюда можно получить (метод простой итерации) численный способ (Шелейховского--Брэгмана--Синхорна) поиска матрицы корреспонденций
\begin{align*}
    d^{rtk}_{ij}\simeq [A^{rt}_{i}]^N [B_j^r]^N Q^{rt}_i D_j^r \exp\left(-\beta^r T^{tk}_{ij}\right),\\
    [A^{rt}_i]^{t'+1} = \left(\sum_{k\in Z(t);~j} [B_j^r]^{t'} D_j^r\exp\left(-\beta^r T^{tk}_{ij}\right)\right)^{-1},\\
    [B_j^r]^{t'+1} = \left(\sum_{t\in M(r);~k\in Z(t);~j} [A^{rt}_i]^{t'+1} Q^{rt}_i\exp\left(-\beta^r T^{tk}_{ij}\right)\right)^{-1}.
\end{align*}
Заметим также, что из отмеченного представления сразу следует, что 
$$\frac{d^{rtk}_{ij}}{\sum_{k'\in Z(t)} d^{rtk'}_{ij}} = \frac{\exp\left(-\beta^r T^{tk}_{ij}\right)}{\sum_{k'\in Z(t)}\exp\left(-\beta^r T^{tk'}_{ij}\right)}.$$
Приведенная формула означает, что можно считать $|Z(t)| = 1$ и просто использовать после блока расчета корреспонденций логит-модель для расщепления пассажиров по типу \cite{Sandholm2010}. 
\fi

\section{Модели равновесного распределения транспортных потоков по путям}\label{section_3}
Матрица корреспонденций $\{d^{rt}_{ij}\}_{(i,j)\in OD}$ порождает (вообще говоря, неоднозначно) вектор распределения потоков по путям $x^{rt}$. Неоднозначность заключается в том, что балансовые ограничения, которые возникают на $x^{rt} \in X^{rt}(d)$:
\begin{equation*}\label{Xd}
  x^{rtk}\ge 0:\quad \forall (i,j)\in OD \to  \sum_{k\in M(r)}\sum_{p\in P^{b(k)}_{ij}} x^{rtk}_p = d^{rt}_{ij},
\end{equation*}
как правило, не определяют вектор $x^{rt}$ однозначно. Кроме того, размерность векторов  $x^{rtk}$ может быть очень большой (число путей может быть очень большим), поэтому удобно будет сформулировать модель равновесного распределения потоков по путям в двойственной форме, подобно \cite{gasnikov2020posobie}. Но сначала выпишем прямую задачу поиска равновесия в модели равновесного распределения потоков по путям, которая получается исходя из потенциальности соответствующей игры загрузки. В такой популяционной <<игре>> тип игроков определяется типом корреспонденции $(i,j)$, которой они принадлежат, слоем спроса $r$ и типом игрока $t$, а стратегией игрока является выбор транспортного средства $k\in Z(t)$ и маршрута $p\in P^{b(k)}_{ij}$ \cite{Sandholm2010,gasnikov2020posobie,mazalov2022}.

\iffalse
Для этого сделаем дополнительное предположение, см. \cite{Shvetsov2003}, заключающееся в том, что (нумерация соответствует примеру выше)
\begin{center}
$\tau^1_e(f_e):=\tau_e(f_e)+c^1_e$, $\tau^3_e(f_e):=\tau_e(f_e)+c^3_e$. 
\end{center}
Поток на ребре есть сумма потоков по путям, проходящем через заданное ребро. Будем также считать, что для общественного транспорта известны $T^2_{ij}$, например, рассчитаны согласно модели оптимальных стратегий Флориана \cite{Shvetsov2003}. Введем также поправочный коэффициент $\kappa_3$ (для остальных значений $k$ будем считать его равным $1$), отвечающий тому, сколько места на дороге занимает служебный транспорт, отвечающий $k = 3$. Например, если это автобус, то $\kappa_3\simeq 3$.
\fi

Введем функцию (опуская возможные различные индексы и волны) $$\sigma_e(f_e) = \int_0^{f_e} \tau_e(z)dz$$ и сопряженную к ней функцию по Фенхелю--Рокафеллару: $$\sigma_e^*(t_e) = \max_{f_e \ge 0} \left\{t_e f_e - \sigma_e(f_e) \right\}.$$ Это определение делается для того, чтобы <<обратить>> формулу $t_e = \tau_e(f_e)$. А именно, $t_e = \tau_e(f_e)$ тогда и только тогда, когда $f_e = d\sigma_{e}^{*}(t_e)/dt_e$.

Введем функционал
$$\Psi(f) = \sum_{b\in B;~e \in E} \sigma_e(f_{e,b}) + \sum_{k\in K;~e \in E} \tilde{\sigma}_e^k(\tilde{f}_{e}^k).$$

Несложно проверить, что этот функционал будет потенциалом в описанной выше популяционной игре. Действительно, учитывая \eqref{t(f)}, \eqref{f(x)}, получаем свойство потенциальности
$$\frac{\partial \Psi (f(x))}{\partial x^{rtk}_p} = \sum_{e\in E} \delta_{ep}t_e^k(f(x)),$$
поскольку правая часть, как раз есть затраты на пути $p$.

Принцип Нэша--Вардропа (что каждый игрок выбирает наилучшую для себя стратегию / маршрут при заданных стратегиях остальных игроков) в данном случае можно записать как  задачу оптимизации 
\begin{equation}\label{primal}
    \min_{x^{rt}\in X^{rt}(d), r, t\in M(r)} \Psi(f(x)).
\end{equation}
 Более того, эта задача будет выпуклой, если функции $\tau_e(f_e)$ не убывают, что является естественным предположением, см., например, \eqref{BPR}. Поэтому, вместо решения исходной задачи, можно перейти к решению двойственной, используя обозначение  \eqref{T}.
 \begin{theorem}\label{Th:dual}
 Двойственная задач для задачи \eqref{primal} будет иметь вид
 \begin{equation}\label{dual}
    \max_{t=\left\{t_{b(k)},\tilde{t}^k\right\}_{k \in K}\ge 0}\left\{\sum_{r;~t\in M(r);~(i,j)\in OD} d_{ij}^{rt} \bar{T}_{ij}^t -\sum_{b\in B;~e \in E} \sigma_e^*(t_{e,b(k)}) - \sum_{k\in K;~e \in E} \tilde{\sigma}_e^{*,k}(\tilde{t}_{e}^k)\right\},
 \end{equation}
где
\begin{equation}\label{min}
\bar{T}_{ij}^t = \min\left\{ \left\{T_{ij}^k(t_{b(k)}+\tilde{t}^k)\right\}_{k\in Z(t)} \right\}.    
\end{equation}
Причем, (формула Демьянова--Данскина) 
\demyan{$$f_e^k = \frac{\partial}{\partial \tilde{t}^k_e} \left\{\sum_{r;~t\in M(r);~(i,j)\in OD} d_{ij}^{rt} \min\left\{ \left\{T_{ij}^k(t_{b(k)}+\tilde{t}^k)\right\}_{k\in Z(t)} \right\}\right\}_{t=t_*}$$}
\demyan{где}  
% $\Phi (t)$ -- двойственный функционал в задаче \eqref{dual}, а 
$t_*$ -- решение задачи \eqref{dual}. Более строго, это равенство следует писать в субдифференциальной форме (детали см., например, \cite{gasnikov2020posobie}).
 \end{theorem}
 Вместо равновесия Нэша--Вардропа можно рассматривать стохастические равновесия, в основе которых лежит выбор игроками не наилучшей стратегии, а рандомизированный выбор стратегии согласно логит-распределению. В таком случае $\min$ можно заменить на его сглаженную версию $\text{softmin}$ \cite{gasnikov2020posobie}. Аналогичный шаг можно сделать и в определении $T^k_{ij}(t^k)$ \eqref{T}. 

 Отметим также, что если $\tilde{\tau}^k_e(f) \equiv c_e^k $, то можно считать $\tilde{t}_e^k = c_e^k$ и исключить соответствующую переменную из двойственной задачи.

 Важным наблюдением является возможность предельного перехода (предел <<стабильной динамики>>) $\mu\to 0+$ и \eqref{BPR} в задачах \eqref{primal} и \eqref{dual}. Такой предельный переход впервые был предложен для более частной модели в работе \cite{gasnikov2014matmod}. Более современное изложение можно найти в работе \cite{kotlyarova2021}.

 \section{Двухстадийная модель}
Стандартный способ поиска равновесий в многостадийных транспортных моделях предполагает последовательную прогонку (отрешивание) двух блоков (двух задач) \eqref{Wilson} и \eqref{primal}. Из решения \eqref{Wilson} находим зависимость $d(T)$, а из решения \eqref{primal} находим зависимость $T(t(f(x(d))))$. Неподвижная точка такой прогонки и будет искомым равновесием. На практике именно такая процедура обычно и реализуется \cite{Ortuzar2002,gasnikov2020posobie,Boyles2020}. В частности, насколько нам известно, в России все основные расчеты в самых разных организация (Генпланах, Департаментах транспорта и т.п.) проводились на основе пакетов транспортного моделирования, в основу которых была положена именно такая прогонка. Однако, чтобы такая процедура сходилась на практике часто необходимо достаточно удачно выбрать точку старта. Более надежный численный способ поиска равновесия в двухстадийной модели заключается в том, чтобы посмотреть на задачи  \eqref{Wilson} и \eqref{dual}, и попытаться объединить эти две задачи оптимизации в одну седловую задачу, учитывая их структуры:
$$\min_{d\in(l,w)} G(d,\bar{T}(t)) + g(d),$$
$$\max_{t\ge 0} \hat{G}(d,\bar{T}(t)) - h(t).$$
Если $G=\hat{G}$ (для этого нужно, чтобы все $\beta_r = 1$, на самом деле, не обязательно, чтобы именно 1, но важно, чтобы были равны между собой), то совместное решение этих двух задач можно найти из решения седловой (выпукло-вогнутой) задачи
$$\min_{d\in(l,w)}\max_{t\ge 0} G(d,\bar{T}(t)) + g(d) - h(t),$$
которую, в свою очередь, можно переписать как (теорема фон Неймана--Сиона--Какутани)
$$\max_{t\ge 0}\min_{d\in(l,w)} G(d,\bar{T}(t)) + g(d) - h(t).$$
И уже для последней седловой задачи можно построить двойственную по части переменных $d$. В результате получается задача (вогнутой) оптимизации с двумя блоками переменных: $t$ и блок двойственных переменных для $d$ (множители Лагранжа к ограничениям $d \in (l,w)$). Мы не будем здесь приводить детали, отметим лишь, что все описанные выкладки будут близки к соответствующим выкладкам работы \cite{gasnikov2020posobie}, в которой рассматривался частный случай описанной модели.

По-видимому, впервые такой способ (для более простой модели) был обнаружен в работе \cite{Evans1976} (более современное изложение см. в \cite{DeCea2005}) и независимо переоткрыт в большей общности и с эволюционной интерпретацией в работе \cite{gasnikov2014matmod} (более современное изложение см. в \cite{gasnikov2020posobie}).  

В настоящей работе также были произведены различные численные эксперименты для графа г. Москвы и области до бетонного кольца (данные предоставили коллеги из Российского университета транспорта) с числом ребер $m\simeq 10^5$. Использованный для решения полученной двойственной задачи прямо-двойственный универсальный ускоренный метод Нестерова (в версии из статьи \cite{gasnikov2018jvm}) в сочетании с методом Синхорна--Брэгмана--Шелейховского \cite{Peyre2019} (метод балансировки) для внутренней задачи (по двойственным к $d \in (l,w)$) переменным) показал себя наилучшим образом, и уменьшил зазор двойственности на 6 порядков за 2 часа работы на современном многоядерном ноутбуке  \cite{kubentaeva2022code}.
% \todo[inline]{У нас GPU не используется}

В заключение заметим, что из-за структуры итоговой двойственной задачи, возможно, ускоренные универсальные блочно-покомпонентные методы \cite{gasnikov2016coordinate} проявят себя еще лучше, чем просто ускоренный универсальный метод в сочетании с методом балансировки. В будущем планируется проверить эту гипотезу.

%\end{fulltext}

%=================Список литературы====================

Исследование в части теории было проведено в рамках выполнения государственного задания Министерства науки и высшего образования Российской Федерации (проект No 0714-2020-0005) совместно с Российским университетом транспорта, любезно предоставившим авторскому коллективу данные для расчетов. В частности, персональная благодарность Алексею Вячеславовичу Шурупову, Владимиру Ивановичу Швецову и Леониду Михайловичу Барышеву, а также Владимиру Викторовичу Мазалову за ценный совет, который помог лучше оформить результат \eqref{primal}. 

Практическая часть исследований была выполнена при поддержке ежегодного дохода ФЦК МФТИ (целевого капитала № 5 на развитие направлений искусственного интеллекта и машинного обучения в МФТИ).

Гасникова Евгения Владимировна, Московский Физико-Технический Институт, 141701, Московская область, г. Долгопрудный, Институтский переулок, д.9. Тел. +7 (926) 855 70 55, e-mail: egasnikov@yandex.ru.

Гасников Александр Владимирович, Московский Физико-Технический Институт, 141701, Московская область, г. Долгопрудный, Институтский переулок, д.9.; Институт Проблем Передачи Информации РАН, 127994, г. Москва, ГСП-4, Большой Каретный переулок, 19, стр. 1; Кавказский математический центр Адыгейского государственного университета, 385000, Республика Адыгея, г. Майкоп, ул. Первомайская, д. 208. Тел. +7 (905) 780 69 74, e-mail: gasnikov@yandex.ru.

Ярмошик Демьян Валерьевич, Московский Физико-Технический Институт, 141701, Московская область, г. Долгопрудный, Институтский переулок, д.9. Тел. +7 (977) 812 13 37, e-mail: yarmoshik.dv@phystech.edu.

Кубентаева Меруза Болатбековна, Московский Физико-Технический Институт, 141701, Московская область, г. Долгопрудный, Институтский переулок, д.9. E-mail: kubentay@gmail.com.

Персиянов Михаил Игоревич, Московский Физико-Технический Институт, 141701, Московская область, г. Долгопрудный, Институтский переулок, д.9.; Институт Проблем Передачи Информации РАН, 127994, г. Москва, ГСП-4, Большой Каретный переулок, 19, стр. 1; Тел. +7 (921) 566 39 31, e-mail: persiianov.mi@phystech.edu.

Подлипнова Ирина Вячеславовна, Московский Физико-Технический Институт, 141701, Московская область, г. Долгопрудный, Институтский переулок, д.9. Тел. +7 (937) 070 85 70, e-mail: podlipnova.iv@phystech.edu. 

Подобная Елена Дмитриевна, Московский Физико-Технический Институт, 141701, Московская область, г. Долгопрудный, Институтский переулок, д.9. Тел. +7 (916) 574-43-42, e-mail: epodobnaya@gmail.com.

Матюхин Владислав Вячеславович, Московский Физико-Технический Институт, 141701, Московская область, г. Долгопрудный, Институтский переулок, д.9. Тел. +7 (925) 991 59 49, e-mail: matyukhin@phystech.edu.


\newpage

\begin{thebibliography}{99}
\bibitem%[Вильсон 1978]
{Wilson1978}
	\textit{Вильсон~А.\,Дж.} Энтропийные методы моделирования сложных систем~// М.: Наука,~--- 1978. 
% {\footnotesize{\it Wilson~A.\,G.} 
% Entropiinye metody modelirovaniya slozhnykh sistem 
% [Entropy in urban and regional modeling]~// M.: Nauka, 1978  (in Russian). \par}
%~// Routledge,~--- 2011. ).\par}
%done ->
% \bibitem[Бабичева и др., 2015]{babicheva2015mipt}
% 	\textit{Бабичева Т.\,С. и др.} Двухстадийная модель равновесного распределения транспортных потоков~// Труды Московского физико-технического института. – 2015. – Т. 7. – №. 3 (27) - С 31-34. \\
% {\footnotesize{\it Babicheva T.\,S.} 
% Dvuhstadijnaya model' ravnovesnogo raspredeleniya transportnyh potokov 
% [Two-stage model of equilibrium distributions of traffic flows]~// 
% Trudy Moskovskogo fiziko-tekhnicheskogo instituta. – 2015. – Vol. 7. – №. 3 (27) - P. 31-34  (in Russian). \par}
%done ->
% \bibitem[Баймурзина и др., 2019]{baimurzina2019jvm}
% 	\textit{Баймурзина Д.\,Р. и др.} Универсальный метод поиска равновесий и стохастических равновесий в транспортных сетях~// Журнал вычислительной математики и математической физики. – 2019. – Т. 59. – №. 1. – С. 21-36. \\
% {\footnotesize{ \it Bajmurzina D.\,R.}
% Universal'nyj metod poiska ravnovesij i stohasticheskih ravnovesij v transportnyh setyah 
% [Universal method of searching for equilibria and stochastic equilibria in transportation networks]~// Computational Mathematics and Mathematical Physics.  – 2019. – Vol. 59. – №. 1. – P. 21-36 (in Russian). \par}
% %done ->	
% \bibitem[Гасников, 2015]{gasnikov2015matmod}
% 	\textit{Гасников А.\,В.} Об эффективной вычислимости конкурентных равновесий в транспортно-экономических моделях //Математическое моделирование. – 2015. – Т. 27. – №. 12. – С. 121-136. \\
% 	{\footnotesize{\it Gasnikov A.\,V.} 
% 	Ob effektivnoj vychislimosti konkurentnyh ravnovesij v transportno-ekonomicheskih modelyah. 
% 	[Reduction of searching competetive equillibrium to the minimax problem in application to different network problems]~//
% 	Mathematical Models and Computer Simulations. – 2015. – V. 27. – №. 12. – P. 121-136 (in Russian). \par}
%done ->	
% \bibitem[Гасников, 2016]{gasnikov2016diss}
% 	\textit{Гасников~А.\,В.} Эффективные численные методы поиска равновесий в больших транспортных сетях.~//Диссертация на соискание степени д.ф.-м.н. по специальности 05.13.18–Математическое моделирование, численные методы, комплексы программ. ~--- М.: МФТИ, ~--- 2016. ~--- 487 с. \\
% 	{\footnotesize{\it Gasnikov~A.\,V.} 
% 	Effektivnye chislennye metody poiska ravnovesii v bol'shikh transportnykh setyakh. 
% 	[Efficient numerical methods for searching equillibriums in large transport networks]~// M: MFTI,~--- 2016.~--- 487 p (in Russian). \par}
%done ->
% \bibitem[Гасников, 2021]{gasnikov2021book}
% 	\textit{Гасников~А.\,В.} Современные численные методы оптимизации. Метод универсального градиентного спуска. //  М.: МЦНМО,~--- 2021.
% 	{\footnotesize{\it Gasnikov~A.\,V.} Sovremennye chislennye metody optimizatsii. Metod universal'nogo gradientnogo spuska. [Modern numerical optimization methods. The universal gradient descent method.] // M: MCCME,~--- 2021(in Russian) \par}
%done ->	
\bibitem%[Гасников и др., 2013]
{gasnikov2013book}
	\textit{Гасников~А.\,В., Кленов~С.\,Л., Нурминский~Е.\,А., Холодов~Я.\,А., Шамрай~Н.\,Б.} Введение в математическое моделирование транспортных потоков. Под ред. А.В. Гасникова с приложениями М.Л. Бланка, К.В. Воронцова и Ю.В. Чеховича, Е.В. Гасниковой, А.А. Замятина и В.А. Малышева, А.В. Колесникова, Ю.Е. Нестерова и С.В. Шпирко, А.М. Райгородского, с предисловием руководителя департамента транспорта г. Москвы М.С. Ликсутова.~// М.: МЦНМО, -- 2013. -- 427 стр., 2-е изд. 
 % \\
	% {\footnotesize{\it Gasnikov~A.\,V. et al.} 
	% Vvedenie v matematicheskoe modelirovanie transportnykh potokov 
	% [Introduction to the mathematical modeling of traffic flows].  Eds. A.V. Gasnikov.// Moscow:MCCME, -- 2013 (in Russian). \par}
%done ->	
\bibitem%[Гасников и др., 2014]
{gasnikov2014matmod}
	\textit{Гасников А.\,В. и др.} О трехстадийной версии модели стационарной динамики транспортных потоков~// Математическое моделирование. – 2014. – Т. 26. – №. 6. – С. 34-70.
	% {\footnotesize{\it Gasnikov~A.\,V. et al.} 
	% O trekhstadijnoj versii modeli stacionarnoj dinamiki transportnyh potokov
	% [On the three-stage version of stable dynamic model]~// 
	% Mathematical Models and Computer Simulations – 2014. – V. 26. – №. 6. – P. 34-70 (in Russian). \par}
%done ->	
% \bibitem[Гасников и др., 2015]{gasnikov2015MIPTershov}
% 	\textit{Гасников А.\,В. и др.} Поиск стохастических равновесий в транспортных моделях равновесного распределения потоков~// Труды Московского физико-технического института. – 2015. – Т. 7. – №. 4 (28) - С. 143-155. \\ 
% 	{\footnotesize{\it Gasnikov~A.\,V. et al.} 
% 	Poisk stohasticheskih ravnovesij v transportnyh modelyah ravnovesnogo raspredeleniya potokov
% 	[Search for the stochastic equilibria in the transport models of equilibrium flow distribution]~//
% 	Trudy Moskovskogo fiziko-tekhnicheskogo instituta. – 2015. – V. 7. – №. 4 (28) - P. 143-155 (in Russian). \par}

%done ->	
\bibitem%[Гасников и др., 2016]
{gasnikov2016} 
    \textit{Гасников~А.\,В. и др.} Эволюционные выводы энтропийной модели расчета матрицы корреспонденций // Математическое моделирование~--- 2016.~--- Т. 28, № 4~--- С. 111-124. 
    % \\
    % {\footnotesize{\it Gasnikov, A. \, V., Gasnikova, E. \,V., Mendel’, M. \,A., Chepurchenko, K. \,V.}
    % Evolyutsionnye vyvody entropiinoi modeli rascheta matritsy korrespondentsii 
    % [Evolutionary interpretations of entropy model for correspondence matrix calculation]~// 
    % Mathematical Models and Computer Simulations~--- 2016.~--- Vol. 28, No. 4~--- P. 111-124 (in Russian).\par}
%done ->    
% \bibitem[Гасников и др., 2021]{gasnikov2021} 
%     \textit{Гасников~А.\,В. и др.} Ускоренный метаалгоритм для задач выпуклой оптимизации~// Журнал вычислительной математики и математической физики~--- 2021.~--- Т. 61, № 1. \\
%     {\footnotesize{\it Gasnikov~A.\,V. et al.} 
%     Uskorennyj metaalgoritm dlya zadach vypukloj optimizacii 
%     [Accelerated meta-algorithm for convex optimization]~// 
%     Computational Mathematics and Mathematical Physics ~--- 2021.~--- T. 61, № 1  (in Russian). \par}
%done ->	
% \bibitem[Гасников--Гасникова, 2010]{gasnikov2010}
% 	\textit{Гасников А.\,В., Гасникова Е.\,В.} О возможной динамике в модели расчета матрицы корреспонденций (А. Дж. Вильсона)~// Труды Московского физико-технического института. – 2010. – Т. 2. – №. 4~--- С. 45-52. \\
% 	{\footnotesize{\it Gasnikov~A.\,V., Gasnikova E.\,V.} 
% 	O vozmozhnoj dinamike v modeli rascheta matricy korrespondencij (A. Dzh. Vil'sona)
% 	[On the possible dynamics in the model for calculating the matrix of correspondences (A.J. Wilson)]~//
% 	Trudy Moskovskogo fiziko-tekhnicheskogo instituta. – 2010. – V. 2. – №. 4~--- P. 45-52  (in Russian). \par}
% %done ->	
% \bibitem[Гасников--Гасникова и др., 2015]{gasnikov2015usik}
% 	\textit{Гасников А. В. и др.} О связи моделей дискретного выбора с разномасштабными по времени популяционными играми загрузок~// Труды Московского физико-технического института. – 2015. – Т. 7. – №. 4 (28)~--- С. 129-142.
% 	{\footnotesize{\it Gasnikov~A.\,V., Gasnikova E.\,V., et al.} 
% 	O svyazi modelej diskretnogo vybora s raznomasshtabnymi po vremeni populyacionnymi igrami zagruzok
% 	[Searching of equilibriums in hierarchical congestion population games]~//
% 	Trudy Moskovskogo fiziko-tekhnicheskogo instituta. – 2015. – V. 7. – №. 4 (28)~--- P. 129-142 (in Russian). \par}
% %done ->	
\bibitem%[Гасников--Гасникова, 2020]
{gasnikov2020posobie}
	\textit{Гасников~А.\,В., Гасникова ~Е.\,В.} Модели равновесного распределения потоков в больших сетях~// М.: МФТИ, ~--- 2020. ~--- 204 с. 
 % \\
	% {\footnotesize{\it Gasnikov~A.\,V., Gasnikova E.\,V.} 
	% Modeli ravnovesnogo raspredeleniya potokov v bol'shikh setyakh
	% [Traffic assignment models. Numerical aspects]~//
	% M: MFTI,~--- 2020.~--- 204 p  (in Russian). \par}
%done ->	
\bibitem%[Гасников и др., 2016]
{gasnikov2016coordinate}
\textit{Гасников А. В., Двуреченский П. Е., Усманова И. Н.} О нетривиальности быстрых (ускоренных) рандомизированных методов // Труды Московского физико-технического института. – 2016. – Т. 8. – №. 2 (30). – С. 67-100.

% \bibitem[Гасников--Дорн и др., 2016]{gasnikov2016matmod}
% 	\textit{Гасников А.\,В. и др.} Численные методы поиска равновесного распределения потоков в модели Бэкмана и в модели стабильной динамики~// Математическое моделирование. – 2016. – Т. 28. – №. 10. – С. 40-64. \\
% 	{\footnotesize{\it Gasnikov~A.\,V. et al.} 
% 	Chislennye metody poiska ravnovesnogo raspredeleniya potokov v modeli Bekmana i v modeli stabil'noj dinamiki 
% 	[Numerical methods for the problem of traffic flow equilibrium in the Beckmann and the stable dynamic models]~//
% 	Mathematical Models and Computer Simulations.  – 2016. – V. 28. – №. 10. – P. 40-64  (in Russian). \par}
%done ->	
% \bibitem[Гасников--Кубентаева, 2018]{gasnikov2018kim}
% 	\textit{Гасников А.\,В., Кубентаева М.\,Б.} Поиск стохастических равновесий в транспортных сетях с помощью универсального прямо-двойственного градиентного метода~//Компьютерные исследования и моделирование. – 2018. – Т. 10. – №. 3. – С. 335-345. \\
% 	{\footnotesize{\it Gasnikov~A.\,V., Kubentaeva M.\,B.} 
% 	Poisk stohasticheskih ravnovesij v transportnyh setyah s pomoshch'yu universal'nogo pryamo-dvojstvennogo gradientnogo metoda 
% 	[Searching stochastic equilibria in transport networks by universal primal-dual gradient method]~//
% 	Computer Research and Modeling. – 2018. – V. 10. – №. 3. – P. 335-345 (in Russian). \par}	
% %done ->	
% \bibitem[Гасников--Нестеров и др., 2015]{gasnikov2015MIPTnesterov}
% 	\textit{Гасников А.\,В. и др.} Поиск равновесий в многостадийных транспортных моделях~//
% 	Труды Московского физико-технического института. – 2015. – Т. 7. – №. 4 (28) - С. 143-155. \\
% 	{\footnotesize{\it Gasnikov~A.\,V. et al.}
% 	Poisk ravnovesij v mnogostadijnyh transportnyh modelyah
% 	[Searching equilibria in multi-stage transport models]~// 
% 	Trudy Moskovskogo fiziko-tekhnicheskogo instituta. – 2015. – V. 7. – №. 4 (28) - P. 143-155 (in Russian). \par}	
%done ->	
% \bibitem[Гасников--Нестеров и др., 2016]{gasnikov2016jvm}
% 	\textit{Гасников А.\,В. и др.} Об эффективных численных методах решения задач энтропийно-линейного программирования~//
% 	Журнал вычислительной математики и математической физики. – 2016. – Т. 56. – №. 4. – С. 523-534. \\
% 	{\footnotesize
% 	{\it Gasnikov~A.\,V. et al.} 
% 	Ob effektivnyh chislennyh metodah resheniya zadach entropijno-linejnogo programmirovaniya 
% 	[Efficient numerical methods for entropy-linear programming problems]~//
% 	Computational Mathematics and Mathematical Physics. – 2016. – T. 56. – №. 4. – S. 523-534 (in Russian). \par}
\bibitem%[Гасников--Нестеров, 2018]
{gasnikov2018jvm}
	\textit{Гасников А.\,В., Нестеров Ю.\,Е.} Универсальный метод для задач стохастической композитной оптимизации~// 
	Журнал вычислительной математики и математической физики. – 2018. – Т. 58. – №. 1. – С. 51-68. \\
	% {\footnotesize
	% {\it Gasnikov~A.\,V., Nesterov Y.\,E }  
	% Universal'nyj metod dlya zadach stohasticheskoj kompozitnoj optimizacii
	% [Universal method for stochastic composite optimization problems]~//
	% Computational Mathematics and Mathematical Physics. – 2018. – V. 58. – №. 1. – P. 51-68 (in Russian). \par}
% \bibitem[ДаДата]{Dadata2020data}
% 	\textit{ДаДата.} {\footnotesize Dadata (in Russian) } 
% 	[Electronic resource]: \url{https://dadata.ru/api/geocode/} (accessed 16.02.2021).
\bibitem%[Иванова, 2020]
{ivanova2020}
	\textit{Иванова A.\,С. и др.} Калибровка параметров модели расчета матрицы корреспонденций для г. Москвы~//
	COMPUTER. – 2020. – Т. 12. – №. 5. – С. 961-978. 
 % \\
	% {\footnotesize
	% {\it Ivanova A.\,C., et al.} 
	% Ivanova A.\,C., et al. Kalibrovka parametrov modeli rascheta matricy korrespondencij dlya g. Moskvy
	% [Calibration of model parameters for calculating correspondence matrix for Moscow]~//
	% COMPUTER. – 2020. – V. 12. – №. 5. – P. 961-978 (in Russian). \par}	
\bibitem%[Кубентаева, 2022]
{kubentaeva2022code} 
	\textit{Кубентаева М.\,Б.} 
	% {\footnotesize \textit{Kubentaeva M.\,B.} }
	[Electronic resource]: \url{https://github.com/MeruzaKub/TransportNet}.
 
\bibitem%[Мазалов--Чиркова, 2022]
{mazalov2022}
\textit{Мазалов В.В., Чиркова Ю.В.}  Сетевые игры. Лань, 2022.

\bibitem%[Котлярова и др., 2021]
{kotlyarova2021} \textit{Котлярова Е. В. и др.} Поиск равновесий в двухстадийных моделях распределения транспортных потоков по сети // Компьютерные исследования и моделирование. – 2021. – Т. 13. – №. 2. – С. 365-379.

\bibitem%[Котлярова и др., 2022]
{kotlyarova2022} 
\textit{Котлярова Е. В.  и др.} Обоснование связи модели Бэкмана с вырождающимися функциями затрат с моделью стабильной динамики // Компьютерные исследования и моделирование. 2022. Т. 14:2. С. 335–342.	
% \bibitem[Котлярова--Ярмошик, 2020]{kоtlyarova2020code} 
% 	\textit{Котлярова Е.\,В., Ярмошик Д.\,В.} 
% 	{\footnotesize \textit{Kotliarova E.\,V., Yarmoshik D.\,V.} }
% 	[Electronic resource]: \url{https://github.com/tamamolis/TransportNet} (accessed 16.02.2021).	
% \bibitem[Министерство труда и социальной политики Приморского края, 2020]{MinTrudPrim2020data}
% Реестр работодателей~// Интерактивный портал Министерства труда и социальной политики Приморского края.
% {\footnotesize Reestr rabotodatelei 
% [Employers register]~//
% Interactive Portal of Ministry of Labour and Social Policy of Primorski Krai
% (in Russian). } 
% [Electronic resource]: \url{https://soctrud.primorsky.ru/employer/} (accessed 16.02.2021).	
\bibitem%[Швецов, 2003]
{Shvetsov2003}	
	\textit{Швецов В. И.} Математическое моделирование транспортных потоков // Автоматика и телемеханика. – 2003. – №. 11. – С. 3-46.	
% \bibitem[Bertsekas, 2009]{Bertsekas2009}
% 	\textit{Bertsekas D.\,P.} Convex optimization theory. – Belmont : Athena Scientific, 2009. 

 \bibitem%[Boyles, 2020]
 {Boyles2020}
 \textit{Boyles S. D., Lownes N. E., Unnikrishnan A.} Transportation network analysis. Volume I: Static and Dynamic Traffic Assignment, 2020.
 
\bibitem%[De Cea, 2005]
{DeCea2005}
 \textit{De Cea J., Fernandez J. E., Dekock V. and Soto A.} Solving network equilibrium problems on
multimodal urban transportation networks with multiple user classes // Transport Reviews. 2005. V. 25(3). P. 293--
317.

 \bibitem%[Evans, 1976]
 {Evans1976}
 \textit{Evans S. P.} Derivation and analysis of some models for combining trip distribution and
assignment // Transportation Research, 1976. V. 10(1). P. 37-–57.	
% \bibitem[Gasnikov et al., 2018]{Gasnikov2018}
% 	\textit{Gasnikov A.\,V., Gasnikova E.\,V., Nesterov Y.\,E.} Dual methods for finding equilibriums in mixed models of flow distribution in large transportation networks~// Computational Mathematics and Mathematical Physics. – 2018. – V. 58. – No. 9. – P. 1395-1403.
\bibitem%[Gasnikova et al., 2023]
{Gasnikova2023}
\textit{Gasnikova E. et al.} An evolutionary view on equilibrium models of transport flows // Mathematics. – 2023. – V. 11. – no. 4. – P. 858.
% \bibitem[Guminov et al., 2019]{Guminov2019}
% 	\textit{Guminov S. et al.} Accelerated alternating minimization~// ICML 2021; arXiv preprint arXiv:1906.03622. – 2019.
	\iffalse
\bibitem[Dvurechensky et al., 2016]{Dvurechensky2016}
	\textit{Dvurechensky P. et al.} Primal-dual method for searching equilibrium in hierarchical congestion population games~// arXiv preprint arXiv:1606.08988. – 2016.
	
\bibitem[Dvurechensky et al., 2018]{Dvurechensky2018}
	\textit{Dvurechensky P., Gasnikov A., Kroshnin A.} Computational optimal transport: Complexity by accelerated gradient descent is better than by Sinkhorn's algorithm~// arXiv preprint arXiv:1802.04367. – 2018.
	
\bibitem[Dvurechensky et al., 2020]{Dvurechensky2020}
	\textit{Dvurechensky P. et al.} A stable alternative to Sinkhorn’s algorithm for regularized optimal transport~// International Conference on Mathematical Optimization Theory and Operations Research. – Springer, Cham, 2020. – P. 406-423.
	
\bibitem[Kamzolov et al., 2020]{Kamzolov2020}
	\textit{Kamzolov D., Dvurechensky P., Gasnikov A.} Universal intermediate gradient method for convex problems with inexact oracle~// Optimization Methods and Software. – 2020. – P. 1-28.
	
\bibitem[Kroshnin et al., 2019]{Kroshnin2019}
	\textit{Kroshnin A. et al.} On the complexity of approximating Wasserstein barycenters~// International conference on machine learning. – PMLR, 2019. – P. 3530-3540.
	
\bibitem[Kubentayeva--Gasnikov, 2020]{Kubentayeva2020}
	\textit{Kubentayeva M., Gasnikov A.} Finding equilibria in the traffic assignment problem with primal-dual gradient methods for Stable Dynamics model and Beckmann model~// arXiv preprint arXiv:2008.02418. – 2020.
	
\bibitem[Nesterov 2015]{Nesterov2015}
	\textit{Nesterov Y.} Universal gradient methods for convex optimization problems~// Mathematical Programming. – 2015. – Т. 152. – №. 1-2. – С. 381-404.
	
\bibitem[Nesterov 2020]{Nesterov2020}
	\textit{Nesterov Y. et al.} Primal–dual accelerated gradient methods with small-dimensional relaxation oracle~// Optimization Methods and Software. – 2020. – P. 1-38.
	\fi
\bibitem%[Nesterov--de Palma 2003]
{Nesterov2003}
	\textit{Nesterov Y., De Palma A.} Stationary dynamic solutions in congested transportation networks: summary and perspectives~// Networks and spatial economics. – 2003. – Т. 3. – No. 3. – P. 371-395.

% \bibitem[Nesterov--Stich, 2017]{Nesterov2017}
% \textit{Nesterov Y., Stich S. U.} Efficiency of the accelerated coordinate descent method on structured optimization problems // SIAM Journal on Optimization. – 2017. – V. 27. – no. 1. – P. 110-123.
\bibitem%[Ortúzar, 2002]
{Ortuzar2002}
	\textit{Ortúzar J.\,D., Willumsen L.\,G.} Modelling transport. John Wiley and Sons~// West Sussex, England. – 2002.
	
\bibitem%[Patriksson, 2015]
{Patriksson2015}
	\textit{Patriksson M.} The traffic assignment problem: models and methods. – Courier Dover Publications, 2015.
	
\bibitem%[Peyré--Cuturi, 2019]
{Peyre2019}
	\textit{Peyré G.,  Cuturi M.} Computational Optimal Transport: With Applications to Data Science~// Foundations and Trends® in Machine Learning. – 2019. – V. 11. – №. 5-6. – P. 355-607.

% \bibitem[Sioux Falls]{SiouxFalls}
% 	\textit{Sioux Falls} 
% 	[Electronic resource]: \url{https://en.wikipedia.org/wiki/Sioux_Falls,_South_Dakota} (accessed 16.02.2021).

\bibitem%[Sandholm, 2010]
{Sandholm2010}
	\textit{Sandholm~W. } Population games and evolutionary dynamics~// MIT press, -- 2010.
	
\bibitem%[Stabler B. et al., 2020]
{Stabler2020}
	\textit{Stabler B., Bar-Gera H., Sall E.} Transportation Networks for Research Core Team. Transportation Networks for Research. Accessed Month, Day, Year. [Electronic resource]: 
	\url{https://github.com/bstabler/TransportationNetworks} (accessed 16.02.2021).
	
% \bibitem[Stonyakin et al., 2020]{Stonyakin2020}
% 	\textit{Stonyakin F.\,S. et al.} Gradient methods for problems with inexact model of the objective~// International Conference on Mathematical Optimization Theory and Operations Research. – Springer, Cham, 2019. – P. 97-114.
	
% \bibitem[Tupitsa et al., 2020]{Tupitsa2020}
% 	\textit{Tupitsa N. et al.} Strongly convex optimization for the dual formulation of optimal transport~// International Conference on Mathematical Optimization Theory and Operations Research. – Springer, Cham, 2020. – P. 192-204.

\end{thebibliography}
\end{document}